\documentclass[10pt]{amsart}

\usepackage{latexsym}
\usepackage{amsmath}
\usepackage{amssymb,amsthm}
\usepackage[english]{babel}
\usepackage{graphicx}
\usepackage{epic, eepic}

\newtheorem{theo}{Theorem}

\newtheorem{lem}{Lemma}
\newtheorem{corol}{Corollary}
\theoremstyle{definition}
\newtheorem*{remarque}{Remark}
\newtheorem*{remarques}{Remarks}

\def\RR{\mathbb{R}}
\def\CC{\mathbb{C}}
\def\NN{\mathbb{N}}
\def\eps{\varepsilon}

\def\supp{{\rm supp}\,}
\newcommand{\trait}[1]{\put(#1,-1){\line(0,1){2}}}
\newcommand{\nm}[2]{\put(#1,-4){\makebox(0,0){#2}}}

\newcommand{\dessineloc}{
\begin{picture}(150,40)(0,-20)
\linethickness{1pt}\thinlines
\put(5,0){\line(1,0){101}}
\put(106,0){\dottedline(0,0)(10,0)}
\put(116,0){\line(1,0){8}}
\put(124,0){\dottedline(0,0)(10,0)}
\put(134,0){\vector(1,0){14}}
\trait{8}
\nm{8}{$\frac{5}{10^{n+3}}$}
\trait{18}
\nm{18}{$\frac{5}{10^{n+2}}$}
\trait{28}
\nm{28}{$\frac{1}{10^{n+1}}$}
\trait{46}
\nm{46}{$\frac{2}{10^{n+1}}$}
\trait{56}
\nm{56}{$\frac{2.5}{10^{n+1}}$}
\trait{66}
\nm{66}{$\frac{3}{10^{n+1}}$}
\trait{76}
\nm{76}{$\frac{3.5}{10^{n+1}}$}
\trait{86}
\nm{86}{$\frac{4}{10^{n+1}}$}
\trait{102}
\nm{102}{$\frac{5}{10^{n+1}}$}
\trait{120}
\nm{120}{$\frac{1}{10^n}$}
\trait{138}
\nm{138}{$\frac{1}{10^{n-1}}$}

\put(150,-1){$r$}
\put(66,4){\makebox(0,0){$\overbrace{\hphantom{\rule{40mm}{2mm}}}$}}
\put(66,8){\makebox(0,0){${\rm supp}\,\chi_n$}}
\put(120,4){\makebox(0,0){$\overbrace{\hphantom{\rule{36mm}{2mm}}}$}}
\put(120,8){\makebox(0,0){${\rm supp}\,\psi_{n-1}$}}
\put(18,4){\makebox(0,0){$\overbrace{\hphantom{\rule{20mm}{2mm}}}$}}
\put(18,8){\makebox(0,0){${\rm supp}\,\psi_{n+1}$}}
\put(66,12){\makebox(0,0){$\overbrace{\hphantom{\rule{20mm}{2mm}}}$}}
\put(66,16){\makebox(0,0){$\chi_n=1$}}
\put(65,-9){\makebox(0,0){$\underbrace{\hphantom{\rule{74mm}{2mm}}}$}}
\put(65,-13){\makebox(0,0){$\psi_n=1$}}
\put(69,-17){\makebox(0,0){$\underbrace{\hphantom{\rule{102mm}{2mm}}}$}}
\put(69,-21){\makebox(0,0){${\rm supp}\,\psi_n$}}
\end{picture}
}
{\newcommand{\dessineW}{
\begin{picture}(140,110)(0,-3)
\linethickness{1pt}
\thinlines
\put(0,0){\vector(1,0){130}}
\put(0,0){\vector(0,1){100}}
\put(132,0){\makebox(0,0){$r$}}
\put(0,103){\makebox(0,0){$W(r)$}}
\spline
(1,99)(2,85)
(3,102)(5,73)
(7,92)(10,55)
(14,80)(19,30)
(28,60)(40,10)
(60,33)(80,0)
(120,5)
\put(20.5,0){\dottedline(0,0)(0,39)}
\put(20.5,-4){\makebox(0,0){$\frac{3}{10^{n+2}}$}}
\put(44,0){\dottedline(0,0)(0,17)}
\put(44,-4){\makebox(0,0){$\frac{3}{10^{n+1}}$}}
\put(0,44){\dottedline(0,0)(23,0)}
\put(-4,44){\makebox(0,0){$\lambda_{n+1}^2$}}
\put(0,22){\dottedline(0,0)(50,0)}
\put(-4,22){\makebox(0,0){$\lambda_{n}^2$}}
\end{picture}

}

\title{\textbf{A Singular Critical Potential For The Schr\"odinger Operator}}
\author{Thomas Duyckaerts}
\address{Universit\'e de Paris-Sud\\ B\^atiment 430\\ 91405 Orsay Cedex} 
\email{Thomas.Duyckaerts@math.u-psud.fr}
\pagestyle{headings}

\begin{document}
\begin{abstract}
We construct a potential $V$ on $\RR^d$, smooth away from one pole, and a sequence of quasi-modes for the operator $-\Delta+V$, which concentrate on this pole. No smoothing effect, Strichartz estimates nor dispersive inequalities hold for the corresponding Schr\"odinger equation.
\end{abstract}
\maketitle
\section{Introduction}

Consider a Schr\"odinger operator:
\begin{equation}
\label{operateur}
P=-\Delta+V,\quad \Delta =\sum_{j=1}^d \frac{\partial^2}{\partial x_j^2}, 
\end{equation}
on $\RR^d$, $d\geq1$, where $V$ is a real potential, small at infinity, and
the initial value problem for the Schr\"odinger equation:
\begin{equation}
\label{LS}
\left\{
\begin{gathered}
i\partial_t U -PU=0\\
U_{\restriction t=0}= U_0 \in L^2(\RR^d).
\end{gathered}
\right.
\end{equation}
For sufficiently smooth potentials $V$, the equation (\ref{LS}) implies a local smoothing effect:
\begin{equation}
\label{regularisant}
||\chi U||_{L^2(]0,T[,H^{1/2}(\RR^d))} \leq C
||U_0||_{L^2(\RR^d)},\; \chi \in C^{\infty}_0(\RR^d),
\end{equation}
where $T$ may be finite, and in some cases infinite, and $C$ may depend
on $V$, $\chi$ and $T$.
This well-known estimate (see \cite{BGT}, \cite{RV} and references therein)
is essential in the study of (\ref{LS}) and related non-linear
equations. In this paper, we investigate the minimal regularity to be assumed on $V$ such that (\ref{regularisant}) is still true. In \cite{RV}, A. Ruiz and
L. Vega consider the equation (\ref{LS}) as a perturbation of the free
Schr\"odinger equation to show an inequality which implies  
(\ref{regularisant}). The following potentials fall within this scope:
\begin{itemize}
\item
$\displaystyle V\in L^{p},\; p\geq d/2.$
\item
Potentials belonging to some Morey-Campanato spaces, (which are larger than $L^{d/2}$) with a smallness condition in these spaces (see \cite{RV} for details). For example, potentials with inverse-square singularities:   
\begin{equation}
\label{potentiel avec pole}
V=\sum_{j=1..N} \frac{a_j}{|x-p_j|^2}+ V_0,\; V_0 \in L^{p},\; p\geq d/2,
\end{equation}
where the $|a_j|$'s are small enough.
\end{itemize}
Let us mention that potentials in the article above may be time-dependent but we shall not develop this point here. 
In \cite{TD}, (see also \cite{BPSTZ}), the author studies potentials of type 
(\ref{potentiel
avec pole}) with regular $V_0$, assuming only the following positivity properties on the $a_j$'s:
$$ a_j+(d/2-1)^2>0,\quad j=1..N,$$
and shows non-trapping type inequality on the truncated resolvent: 
\begin{equation}
\label{resolvante}
 \forall \lambda \in \CC,\; {\rm Re}\, \lambda^2 \geq C,\; {\rm Im}\,\lambda^2
  \neq 0,\quad ||\chi (P-\lambda^2)^{-1} \chi ||_{L^2\rightarrow L^2} \leq \frac{C}{|\lambda|},\quad
\chi  \in C_0^{\infty}(\RR^d),
\end{equation}
One may deduce from (\ref{resolvante}) estimates of type 
(\ref{regularisant}) (cf \cite{BGT}).\par
The preceding works do not study potentials $V\in L^p$, $p<d/2$, and potentials of type (\ref{potentiel avec pole}) but with poles of higher order. 
In this paper, we show that the space $L^{d/2}$ and the inverse-square
singularities are critical, by constructing an unipolar potential whose pole is
of the order of $(\log r)^2/r^2$, $r=|x|$, near $0$, and such that
(\ref{regularisant}), (\ref{resolvante}) and some other usual estimates on
solutions of Schr\"odinger equations do not hold. 
\begin{theo}
Let $d\geq 1$, $N\geq0$ be integers. There exist:
\begin{itemize}
\item a radial, positive potential $V$ on $\RR^d$, which has compact support and such that:
\begin{gather}
\label{HV1}
V \in C^\infty(\RR^d \backslash \{0\})\\
\label{HV2}
\frac{|\log r|^2}{Cr^2}\leq V(r) \leq \frac{C |\log r|^2}{r^2},\; r\leq r_0,
\end{gather}
\item an increasing sequence $(\lambda_n)_{n\geq n_0}$ of positive real numbers, diverging to $+\infty$, 
\item a sequence of radial $C^{\infty}$ functions $(u_n)_{n\geq n_0}$, whose support is of the following form:
\begin{equation*}
\left\{ c_1 \frac{\log \lambda_n}{\lambda_n}\leq r \leq c_2\frac{\log
  \lambda_n}{\lambda_n} \right\}; 
\end{equation*}
\end{itemize}
such that:
\begin{gather}
\label{equation} 
(-\Delta+V)u_n-\lambda_n^2 u_n=f_n \\
\label{thm norme u}
||u_n||_{L^1}=1\\
\label{thm norme f}
\forall j\in \NN,\;||\frac{d^j}{dr^j} f_n||_{L^\infty}=O(\lambda_n^{j-N}),\;
n\rightarrow +\infty.
\end{gather}
\end{theo}

\begin{corol}
Let $N$ be any integer greater than $1$, $P=-\Delta+V$, where $V$ is the
potential of the preceding theorem, and $\chi$ a function in $C_0^{\infty}(\RR^d)$ which does not vanish in $0$. Then:
\begin{itemize}
\item (lack of local Strichartz estimates) 
\begin{equation}
\label{reg 1}
\forall q,q_0 \in [1,+\infty],\; q>q_0,\;\forall T>0,\; \forall C>0,\;  
\exists U_0\in C^{\infty}_0,\;||\chi U(t)||_{L^1(]0,T[,L^q(\RR^d))} > C||U_0||_{L^{q_0}(\RR^d)};
\end{equation}
\item (lack of local smoothing effect)
\begin{equation}
\label{reg 2}
 \forall \sigma>0\;\forall T>0,\; \forall C>0,\;
\exists U_0\in C^{\infty}_0,\; ||\chi U(t)||_{L^1(]0,T[,H^\sigma(\RR^d))} >C ||U_0||_{L^2(\RR^d)};
\end{equation}
\item (lack of local dispersion)
\begin{equation}
\label{dispersion}
\forall q\in ]1,+\infty],\; \forall T>0,\; \forall C>0,\; \exists
U_0\in C^{\infty}_0,\;
||\chi U(T)||_{L^q(\RR^d)} > C||U_0||_{L^{q'}(\RR^d)};
\end{equation}
\item (lack of Strichartz estimates with loss of derivative)
\begin{multline}
\label{Strichartz}
\forall q\in[1,+\infty],\;\forall \sigma\in[0,1],\; \frac 12 -\frac 1q
>\frac{\sigma}{d},\\ \;\forall T>0,\;\forall C>0,\; \exists U_0\in C^{\infty}_0,\;||\chi U(t)||_{L^1(]0,T[,L^q(\RR^d))}> C||U_0||_{D(P^{\sigma/2})}. 
\end{multline}
\end{itemize}
In the preceding statements, we wrote $U(t)$ the solution of the equation (\ref{LS}) with initial value $U_0$ and $q'$ the conjugate exponent of $q$ which is defined by: $1/q+1/q'=1$.
\end{corol}
\begin{remarques}
\begin{itemize}
\item If $d>2$, the hypothesis (\ref{HV1}) and (\ref{HV2}) implie:
$$ V\in \bigcap_{p< d/2} L^p. $$
\item Of course, the sequence $(f_n)_{n\geq n_0}$ invalidates 
(\ref{resolvante})  when $N\geq 2$.
\item It will be clear in the proof of the theorem than one may construct quasi-modes of infinite order (i.e. such that(\ref{thm norme f}) holds with  $O(\lambda_n^{-\infty})$ instead of $O(\lambda_n^{j-N})$) by taking a singularity just a bit stronger: 
$$ \frac{|\log r|^{2+\eps}}{Cr^2}\leq V(r) \leq \frac{C |\log
r|^{2+\eps}}{r^2},\; \eps>0.$$
With a still stronger singularity, one may force $f_n$ to be exponentially decreasing in $-\lambda_n$. 
\item Classical Strichartz estimates:
\begin{gather*}
||u||_{L^p\left(]0,T[,L^q(\RR^d)\right)}\leq C ||u_0||_{L^2(\RR^d)},\;
p>2,\; \frac 2p+ \frac dq=\frac d2,
\end{gather*}
are invalidated by (\ref{reg 1}), with the exception of the case $q=2$, $p=+\infty$ which is always true by the $L^2$ conservation law. Likewise, (\ref{Strichartz}) shows that Strichartz estimates with loss of derivative (see \cite{BGT}) of the form: 
\begin{equation*}
||u||_{L^p\left(]0,T[,L^q(\RR^d)\right)}\leq C
||u_0||_{D(P^{\sigma/2})},\;\frac 12-\frac 1q>\frac \sigma d,\;
\sigma\in [0,1],
\end{equation*}
are false.
In this inequality, the limit case $\sigma/d=1/2-1/q$ is an immediate consequence of the Sobolev inclusion
$$ L^q(\RR^d) \subset H^{\sigma}(\RR^d)\subset D(P^{\sigma/2}),$$
and of the conservation of the norm $D(P^{\sigma/2})$ of any solution of
(\ref{LS}). Notice that this last inclusion implies that (\ref{Strichartz}) is still true in usual Sobolev spaces, i.e when $D(P^{\sigma/2})$ is replaced by $H^{\sigma}(\RR^n)$. 

\item 
The potential $V$ and quasi-modes $u_n$ being of compact support, as close to $0$ as desired, the preceding counter-example is still valid in any regular domain, for a Laplace operator with Dirichlet or Neumann Dirichlet boundary condition.
\item 
Examples of linear Schr\"odinger equations which do not admit the classical
dispersive inequality or Strichartz estimates are given in  \cite{Ban} and
\cite{CZ}. In these two cases, the particular behaviour of the equation
arises from the metric which defines the Laplace operator. Let us mention that the idea to use quasi-modes in relation with Strichartz inequalities is owed to M. Zworski (see \cite{SS}).
\end{itemize}
\end{remarques}
Finally, we would like to point out that for the potential $V$ introduced in the theorem, one may define the operator $P$ without ambiguity. The operator $-\Delta+V$ has a natural meaning on $C^{\infty}_0(\RR^d \backslash\{0\})$, and is essentially self-adjoint on this space, under a positivity asumption implied by (\ref{HV2}). We call $P$ its unique self-adjoint extension, which is of course positive. We refer to Reed and Simon \cite[Th X.30]{RS} for any precision. \par
The author would like to thank Nicolas Burq, its advisor and Clotilde Fermanian Kammerer for their precious help.\par
Section \ref{dem cor} of the paper is devoted to the proof of the corollary, and section \ref{dem theo} to that of the theorem. 

\section{Proof of the corollary.}
\label{dem cor}
According to H\"older's inequality, for any function $v$ on $\RR^d$, with a finite volume support:
\begin{equation}
\label{volume}
||v||_{L^{q_0}}\leq B^{1/q_0-1/q}||v||_{L^q} ,\quad q>q_0,
\end{equation}
where $B$ is the volume of the support.
Let $U_n(t)=e^{-i\lambda_n^2 t} u_n$. Then:
\begin{gather}
\notag
i\partial_t U_n -P U_n  =  -e^{-i\lambda_n^2 t} f_n,\quad
U_{n\restriction t=0} = u_n\\
\label{expression Un}
U_n(t)=e^{-itP}u_n+i e^{-i\lambda_n^2 t} \int_0^t e^{i(s-t)P} f_n ds.
\end{gather}
Assume that (\ref{reg 1}) is false. 
According to (\ref{expression Un}):
$$ \int_0^T ||\chi U_n(t)||_{L^q(\RR^d)}dt\leq \int_0^T ||\chi e^{-itP}u_n||_{L^q(\RR^d)} dt +\int_0^T \int_0^t ||\chi e^{i(s-t)P}f_n||_{L^q(\RR^d)} ds dt.$$
Since the support of $u_n$ concentrates in $0$, and $\chi$ does not vanish
in $0$, the left term of this inequality is, for $n$ large enough, greater than:
\begin{equation*}
\frac 1C \int_0^T ||u_n||_{L^q(\RR^d)} dt\geq \frac 1C \left( \frac{\log \lambda_n}{\lambda_n}\right)^{d(1/q-1/q_0)}||u_n||_{L^{q_0}(\RR^d)}. 
\end{equation*}
Using the negation of (\ref{reg 1}), the right term is dominated by:
\begin{equation*}
||u_n||_{L^{q_0}(\RR^d)}+||f_n||_{L^{q_0}(\RR^d)}.
\end{equation*}
Hence we obtain, using (\ref{thm norme f}):
$$\left( \frac{\lambda_n}{\log \lambda_n} \right)^{d(1/q_0-1/q)}
||u_n||_{L^{q_0}}=O\left(||u_n||_{L^{q_0}}+\lambda_n^{-N}\right),$$
which leads to the announced contradiction since $N\geq 1$ and $q>q_0$, the norm of $u_n$ in $L^{q_0}$ being greater than $1$ by (\ref{thm norme u}).\par
By Sobolev inequalities, (\ref{reg 1}) implies (\ref{reg 2}), taking in (\ref{reg 1}) $q_0=2$ and $q>2$ close enough to $2$. So (\ref{reg 2}) holds.\par
Let us assume that (\ref{dispersion}) is not true. Then we get, by (\ref{expression Un}):
\begin{align*}
||\chi u_n(T)||_{L^q(\RR^d)}&\leq C||u_n||_{L^{q'}(\RR^d)} +\int_0^T
||\chi e^{i(s-T)} f_n||_{L^1(\RR^d)} ds\\
&\leq C\left(||u_n||_{L^{q'}(\RR^d)}+||f_n||_{L^2(\RR^d)}\right).
\end{align*}
To obtain the second inequality we have bounded, up to a multiplicative constant, the $L^{q'}$ norm on the support of $\chi$ by the $L^2$ norm. The contradiction is again a simple consequence of (\ref{volume}) and the fact that $q'<q$.\par
We shall now prove (\ref{Strichartz}). Otherwise, we would have by (\ref{expression Un}):
\begin{equation}
\label{Strichartz un}
\int_0^T ||u_n(t)||_{L^q(\RR^d)} dt \leq
C||u_n||_{D(P^{\sigma/2})}+\int_0^T \int_0^t||e^{i(s-t)P}
f_n||_{D(P^{\sigma/2})} ds\,dt.
\end{equation}
Of course, $e^{i(s-t)P}$ maps isometrically the domain of $P^{\sigma/2}$ onto
itself. We have:
$$ ||f_n||^2_{D(P^{1/2})}=\int |\nabla f_n(x)|^2dx+\int (1+V(x))|f_n(x)|^2 dx.$$
On the support of $f_n$, as on that of $u_n$, $|x|\geq \frac {\log
\lambda_n}{C \lambda_n}$. Using the superior bound of $V$ in (\ref{HV2}), we get, on the support of $f_n$:
$$ |V(r)| \leq C \frac{\lambda_n^2}{(\log \lambda_n)^2} \left(
\log\log\lambda_n-\log \lambda_n\right)^2\leq C\lambda_n^{2}.$$
It follows using (\ref{thm norme f}) and $N\geq 2$ that:
\begin{equation}
\label{majoration fn dans D}
||f_n||_{D(P^{1/2})} \leq C\left(||\nabla
f_n||_{L^2}+\lambda_n||f_n||_{L^2}\right) \leq C\lambda_n^{-1}.
\end{equation} 
Furthermore, the equation $P u_n-\lambda_n^2 u_n=f_n$ implies, using once again the bounds (\ref{thm norme u}) and (\ref{thm norme f}) on the norms of $u_n$ and $f_n$:
\begin{align*}
||u_n||^2_{D(P^{1/2})} & \leq ||f_n||_{L^2} ||u_n||_{L^2}+\lambda_n^2
||u_n||^2_{L^2}\\
& \leq C\lambda_n^2 ||u_n||_{L^2}^2.
\end{align*}
By interpolation on the norms of the left hand side, and since $0\leq\sigma\leq 1$, we get:
\begin{equation}
\label{majoration un dans D}
||u_n||_{D(P^{\sigma/2})}\leq C \lambda_n^\sigma ||u_n||_{L^2}.
\end{equation}
According to (\ref{majoration fn dans D}), (\ref{majoration un dans
D}) and the inequality (\ref{Strichartz un}):
\begin{equation*}
||u_n||_{L^q} \leq C\lambda_n^\sigma ||u_n||_{L^2} +o(1).
\end{equation*}
Hence, with (\ref{volume}):
\begin{equation*}
\left(\frac{\lambda_n}{\log \lambda_n}\right)
^{d/2-d/q} ||u_n||_{L^2}  \leq C \lambda_n^\sigma ||u_n||_{L^2} +o(1),
\end{equation*}
which is absurd since $\frac d2-\frac dq>\sigma$.

\section{D\'emonstration du th\'eor\`eme.}
\label{dem theo}
Denote by $r$ the euclidian norm of $x$ and by $'$ the radial derivative $\frac{d}{dr}$. We would like to find radial functions $V$, $f_n$, $u_n$ such that: 
\begin{equation}
\label{f_n}
 f_n(r)=-u_n''(r)-\frac{d-1}{r} u_n'(r)+V(r) u_n(r) -\lambda_n ^2 u_n(r),
\end{equation}
with $f_n$ small and $\lambda_n$ diverging to infinity. We shall first change functions to get rid of the first-order derivative in this equation. Set:
\begin{equation}
\label{def vn gn W}
u_n=r^{-\frac{d-1}{2}}v_n ,\quad f_n=r^{-\frac{d-1}{2}}g_n,\quad W=V+ \frac{d^2-4d+3}{4 r^2}.
\end{equation}
Thus (\ref{f_n}) becomes:
\begin{equation}
\label{g_n}
\tag{\ref{f_n}'}
g_n(r)=-v_n''+W v_n - \lambda_n^2 v_n.
\end{equation}
Let:
$$y(s)=e^{-\sqrt{s^2+1}},\;
b(s)=-\frac{1}{(s^2+1)^{3/2}}+\frac{s^2}{s^2+1},$$
 which are $C^{\infty}$ solutions of the equation on $\RR$:
\begin{gather}
\label{eqy}
-y''(s)+b(s)y(s)=0\\
\label{inegy}
y(s)>0,\forall j\in \NN,\quad |y^{(j)}(s)|\leq C_je^{-|s|}\\
\label{inegb}
|b(s)|\leq 1.
\end{gather}
We shall write:
\begin{gather*}
y_{\omega,a}(r)=y(\omega(r-a)),\quad b_{\omega,a}(r)=\omega^2 b(\omega(r-a)),
\end{gather*}
where $a$ and $\omega$ are two real parameters. We have:
\begin{equation}
\label{eqy'}
\tag{\ref{eqy}'}
-y''_{\omega,a}+b_{\omega,a} y_{\omega,a}=0;
\end{equation}
Let $q(\lambda)$ be a strictly increasing positive function defined in a neighbourhood of $+\infty$ such that:
\begin{gather}
\label{hypq1}
\lim_{\lambda \rightarrow +\infty} q(\lambda)=+\infty\\
\label{hypq2}
\lim_{\lambda \rightarrow +\infty} \frac{q(\lambda)}{\lambda}=0.
\end{gather}
Let $(\lambda_n)_{n\geq n_0}$ be the sequence defined by the equations:
$$ 10^n=q(\lambda_n),$$ 
so that $\lambda_n$ diverges faster to $+\infty$ than $10^n$.\par
Choose a cutoff function $\chi$:
$
\displaystyle
\chi \in C^{\infty}_0\left(\left]-1,1\right[ 
\right),\; \chi=1 \;{\rm on }\; \left[-\frac 12,\frac 12\right].
$\par
Set:
\begin{gather*}
\psi_n(r)=\chi(10^n r)-\chi(10^{n+1}r),\\
\chi_n(r)=\chi\left( 10^{n+1}\left(r-\frac{3}{10^{n+1}}\right)\right),
\end{gather*}
\setlength{\unitlength}{1mm}
\begin{center}
\begin{figure}
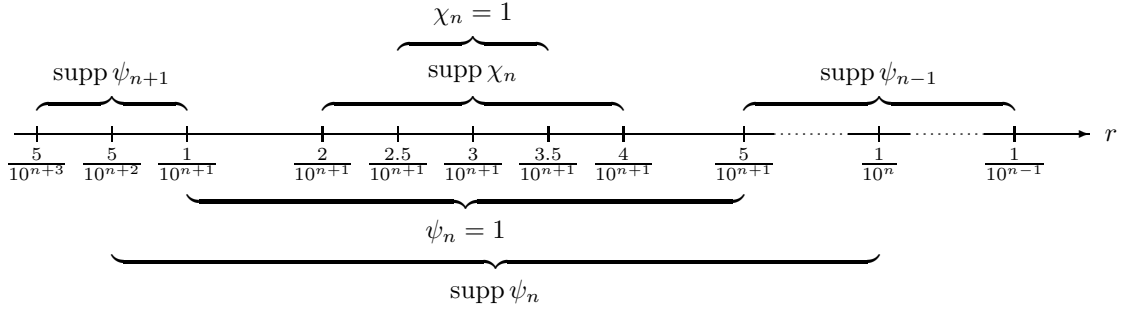

\dessineloc
\caption{Cutoff functions near $3/10^{n+1}$}
\end{figure}
\end{center}
The $\psi_n$'s form a partition of unity near $0$. The support of each $\chi_n$ is included in that of $\psi_n$, away from those of the $\psi_j$'s, $j\neq n$:
\begin{gather}
\notag
\supp \psi_n\subset \left[\frac{5}{10^{n+2}},\frac{1}{10^{n}}\right]\\
\notag
\frac {1}{10^{n+1}}\leq r\leq\frac {5}{10^{n+1}} \Rightarrow \psi_n(r)=1\\
\label{somme psin}
\sum_{n\geq n_0} \psi_n(r)=\chi(10^{n_0}r)\\
\label{support chin}
\supp \chi_n \subset \left[ \frac{2}{10^{n+1}},\frac {4}{10^{n+1}} \right]\subset \{\psi_n=1\}\\
\label{chin=1}
\frac {25}{10^{n+2}}\leq r\leq \frac {35}{10^{n+2}} \Rightarrow \chi_n(r)=1.
\end{gather}
We shall denote by $y_n$ and $b_n$ the following functions:
\begin{gather}
\label{def yn}
 y_n(r)=y_{\frac{\lambda_n}{2},\frac{3}{10^{n+1}}}(r)=y\left( \frac {\lambda_n}{2}\left(r-\frac{3}{10^{n+1}}\right)\right)\\
 b_n(r)=b_{\frac{\lambda_n}{2},\frac{3}{10^{n+1}}}(r)=\frac{\lambda_n^2}{4} b\left(\frac{\lambda_n}{2}\left(r-\frac{3}{10^{n+1}}\right)\right).
\end{gather}
Each of the function $-b_n$ may be seen as a potential well which concentrates (by equation (\ref{eqy'})) the energy of $y_n$ near $3/10^{n+1}$. We shall construct $W$ so that the equation $g_n(r)=0$ is exactly $-v_n''+b_n v_n=0$ on a small interval (including the support of $v_n$) around the point $3/10^{n+1}$. The size of this interval will be of the same order as $10^{-n}$, smaller than the equation parameter $\omega=\lambda_n/2$. We shall choose $v_n$ as a cutoff of $y_n$. The exponential decay of $y$, combined with the difference of order between these two scales will induced small enough error terms for (\ref{thm norme u}) and (\ref{thm norme f}) to hold.
\begin{center}
\begin{figure}
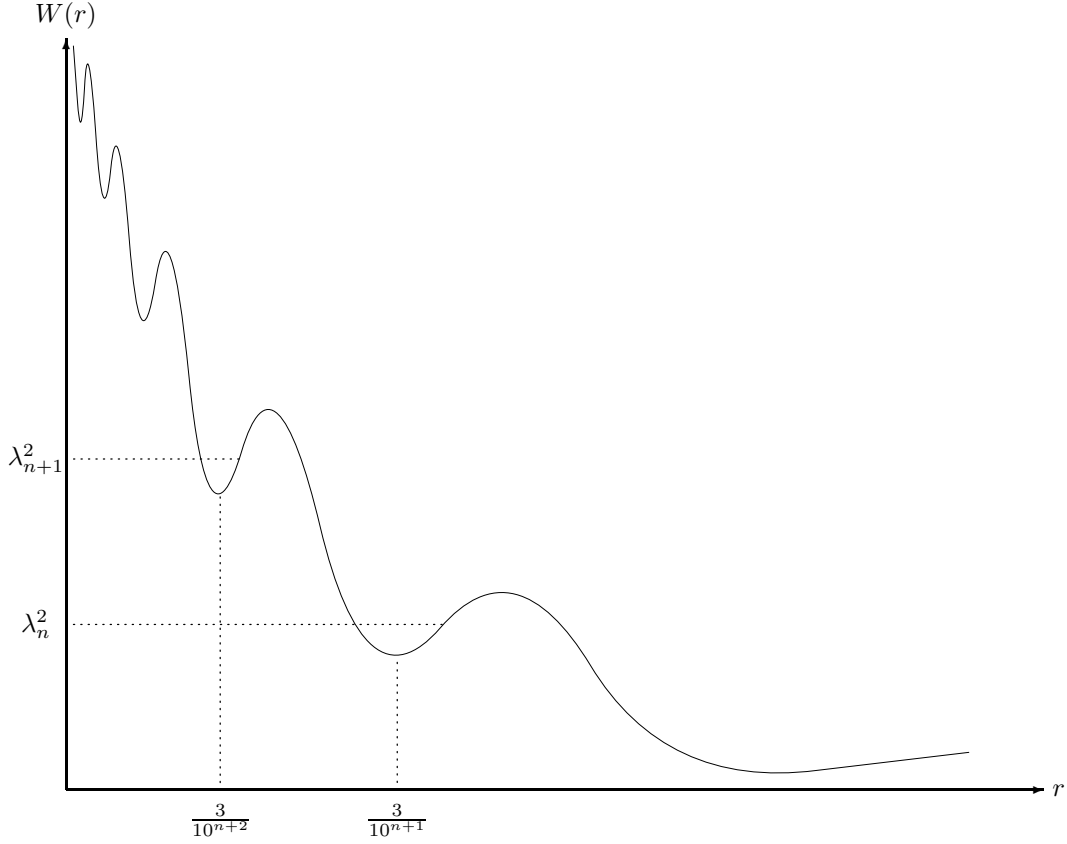

\label{figure 1}
\dessineW
\caption{The potential W}
\end{figure}
\end{center}
Set (see the figure):
\begin{gather}
\label{def W}
W(r)=\sum_{n\geq n_0} \psi_n(r)
\left(b_n(r)+\lambda_n^2\right)\\
v_n(r)= \alpha_n y_n(r)\chi_n(r),
\end{gather}
where the $\alpha_n$'s are strictly positive constants yet to be determined. The functions  $g_n$, $u_n$, $f_n$ and the potential $V$ are thus defined by (\ref{g_n}) and (\ref{def vn gn W}).
\begin{lem}
\label{estimations}
The following inequalities hold:
\begin{gather}
\label{ineg fn}
\forall j\in \NN,\;\exists C>0,\; ||\frac{d^j f_n}{dr^j}||_{L^{\infty}(\RR^d)}\leq C \alpha_n \lambda_n^{1+j}\left(q(\lambda_n)\right)^{\frac{d+1}{2}} 
 e^{-\frac{\lambda_n}{20 q(\lambda_n)}}\\
\label{ineg un}
||u_n||_{L^1(\RR^d) } \geq \frac{\alpha_n}{C} \lambda_n^{-1} \left(q(\lambda_n)\right)^{-\frac{d-1}{2}}.
\end{gather}
\end{lem}
\begin{proof}
According to (\ref{support chin}), on the support of $\chi_n$, $W$ is equal to $b_n+\lambda_n^2$. Hence:
\begin{equation}
\label{expression gn}
g_n(r) =-v_n''(r) + b_n(r) v_n(r).
\end{equation}
Thus, using (\ref{eqy'}):
\begin{equation*}
f_n(r)=r^{-\frac{d-1}{2}} g_n(r)=r^{-\frac{d-1}{2}}\alpha_n \left(y_n(r) \chi_n''(r) +2 y'_n(r) \chi_n'(
r)\right).
\end{equation*}
The derivative of $f_n$ of order $j$ is then of the form:
\begin{equation}
\label{derivee fn}
\frac{d^j f_n}{dr^j}=\sum_{j_1+j_2+j_3=j+1}
\beta_{j_1,j_2,j_3}\frac{d^{j_1}}{dr^{j_1}}y_n 
\frac{d^{j_2}}{dr^{j_2}}\chi_n'\frac{d^{j_3}}{dr^{j_3}} r^{-\frac{d-1}{2}}. 
\end{equation}
On the support of $\chi'_n$:
$$ \left|r- \frac{3}{10^{n+1}}\right| >
\frac{1}{10^{n+1}}=\frac{1}{10q(\lambda_n)}.$$ 
So according to the bounds (\ref{inegy}) of $y$ and its derivatives:
\begin{align*}
|y_n^{(j_1)}(r)|&\leq C\lambda_n^{(j_1)} e^{-\frac{\lambda_n}{2}\left|r-
\frac{3}{10^{n+1}}\right|}\\
&\leq \lambda_n^{j_1} e^{-\frac{\lambda_n}{20q(\lambda_n)}}
\end{align*}
Furthermore,
\begin{alignat*}{3}
|\chi_n^{(j_2+1)}|&\leq C\left(q(\lambda_n)\right)^{j_2+1}&\leq& C\lambda_n^{j_2}q(\lambda_n)\\
\left|\left(\frac{d}{dr}\right)^{j_3}\left(
    r^{-\frac{d-1}{2}}\right)\right|&\leq
C\left(q(\lambda_n)\right)^{\frac{d-1}{2}+j_3}&\leq& C\lambda_n^{j_3}\left(q(\lambda_n)\right)^{\frac{d-1}{2}},
\end{alignat*}
for on the support of $\chi_n$, $r\geq \frac{1}{10^{n+1}}$.
These three inequalities, together with (\ref{derivee fn}), implies
(\ref{ineg fn}).\par
We shall now prove (\ref{ineg un}). By the definition of $y_n$:
$$ \frac{\lambda_n}{2} \left | r-\frac{3}{10^n} \right|\leq \frac 12 \Rightarrow y_n(r)\geq m=\sup_{|s|\leq \frac 12} |y(s)|.$$
Furthermore, if $r$ is as above, and $n$ big enough, then $\chi_n(r)=1$ and so:
$$ u_n(r)=r^{-\frac{d-1}{2}} v_n(r)=r^{-\frac{d-1}{2}} \alpha_n y_n(r)$$
Hence:
\begin{align*}
||u_n||_{L^1} &\geq m \alpha_n \int_{
\lambda_n \left|r-\frac{3}{10^n} \right|\leq 1} r^{-\frac{d-1}{2}} r^{d-1} dr\\
&\geq \frac{\alpha_n}{C} \left(10^{-n}\right)^{\frac{d-1}{2}} \lambda_n^{-1}. 
\end{align*}
\end{proof}
Choose $M>1$ and set: 
\begin{equation}
\label{defq}
q(\lambda)= \frac{\lambda}{M\log \lambda}. 
\end{equation}
The positive function $q$ is strictly increasing for big $\lambda$'s and satisfies (\ref{hypq1}) and (\ref{hypq2}). In addition, lemma \ref{estimations} 
implies, bounding from above $q(\lambda_n)$ by $\lambda_n$:
\begin{gather*}
||\frac{ d^jf_n}{dr^j}||_{L^{\infty}} \leq C \alpha_n \lambda_n^{\frac{d+3}{2}+j-\frac{M}{20}}\\
||u_n||_{L^1} \geq \frac{\alpha_n}{C} \lambda_n^{\frac{d+1}{2}},
\end{gather*}
with new constants $C$, which may depend on $M$.
So the conditions (\ref{thm norme u}) and (\ref{thm norme f}) of the theorem are satisfied if the constants $M$ and $\alpha_n$ are well chosen. The support of $u_n$ is that of $\chi_n$, which is of the desired form:
$$\left \{ c_1\frac{\log \lambda_n}{\lambda_n}\leq r \leq c_2
  \frac{\log\lambda_n}{\lambda_n}\right\},$$
if the support of $\chi$ is taken to be a segment.\par
The assertion (\ref{HV2}) on the potential remains to be checked. We have the following approximation of the inverse function of $q$:
\begin{lem}
\label{equiv q}
Let $q$ be defined by (\ref{defq}). Then:
$$ q(\lambda)\log(q(\lambda)) \sim \frac{\lambda}{M},\; \lambda\rightarrow +\infty.$$
\end{lem}
\begin{proof}
\begin{align*}
q(\lambda)\log q(\lambda)&=q(\lambda) \log\left(\frac{\lambda}{M\log \lambda}\right)\\
&=q(\lambda)\left(\log \lambda-\log\log \lambda -\log M\right)\\
&\sim q(\lambda)\log\lambda,
\end{align*}
when $\lambda$ goes to infinity.
\end{proof}
On the support of $\psi_n$, $q(\lambda_n)=10^n$ and so:
\begin{gather*}
 \frac{1}{20r}\leq q(\lambda_n) \leq \frac{1}{r}\\
\frac{1}{20r}|\log r -\log 20| \leq q(\lambda_n) \log(q(\lambda_n))\leq \frac{1}{r} |\log r|.
\end{gather*}
(We used that $s\mapsto s \log s$ is an increasing function for $s>e^{-1}$).
Using lemma \ref{equiv q} and taking $r$ close enough to $0$ we get: 
$$ C^{-1} \frac 1r |\log r|\leq \lambda_n \leq C \frac{1}{r} |\log r|. $$
Thus, by the definition (\ref{def W}) of $W$, $|b_n|$ being bounded from above by $\lambda_n^2/4$:
$$ C^{-1} \frac {|\log r|^2}{r^2} \sum_{n\geq n_0} \psi_n(r) \leq W(r) \leq C  \frac {|\log r|^2}{r^2} \sum_{n\geq n_0} \psi_n(r), $$
which implies, with (\ref{somme psin}), the inequality (\ref{HV2}) on the potential $W$. Consequently $V$ satisfies the same inequality. 
\begin{remarque} 
To get quasi-modes of infinite order, we would have taken $q(\lambda)=\frac{\lambda}{(\log\lambda)^{1+\eps}}$ and suitably modified lemma \ref{equiv q}.
\end{remarque}

\end{document}


\maketitle
\begin{abstract}
On construit un potentiel $V$ sur $\RR^d$, r\'egulier en dehors d'un p\^ole, et une suite de quasi-modes pour l'op\'erateur $-\Delta+V$, qui se concentre sur ce p\^ole. Les solutions de l'\'equation de Schr\"odinger associ\'ee ne v\'erifient pas d'effet r\'egularisant local, ni d'in\'egalit\'es de Strichartz ou de dispersion. 
\end{abstract}

\section{Introduction}

On consid\`ere un op\'erateur de la forme:
\begin{equation}
\label{operateur}
P=-\Delta+V,\quad \Delta =\sum_{j=1}^d \frac{\partial^2}{\partial x_j^2}, 
\end{equation}
sur $\RR^d$, $d\geq1$, o\`u $V$ est un potentiel r\'eel, petit \`a
l'infini, et l'\'equation de Schr\"odinger:
\begin{equation}
\label{LS}
\left\{
\begin{gathered}
i\partial_t U -PU=0\\
U_{\restriction t=0}= U_0 \in L^2(\RR^d).
\end{gathered}
\right.
\end{equation}
Lorsque $V$ est suffisamment r\'egulier, l'\'equation (\ref{LS})
entra\^ine un effet r\'egularisant local:
\begin{equation}
\label{regularisant}
||\chi u||_{L^2(]0,T[,H^{1/2}(\RR^d))} \leq C
||u_0||_{L^2(\RR^d)},\; \chi \in C^{\infty}_0(\RR^d),
\end{equation}
o\`u selon les cas $T$ est fini ou infini, et $C$ peut d\'ependre
de $V$, $\chi$ et $T$.
Cette estimation classique (cf \cite{BGT},\cite{RV} et les r\'ef\'erences
propos\'ees dans ces deux articles) est fondamentale dans l'\'etude de (\ref{LS})
et des \'equations non-lin\'eaires associ\'ees. 
On s'int\'eresse \`a la r\'egularit\'e minimale de $V$
pour que (\ref{regularisant}) soit encore
v\'erifi\'ee. Dans \cite{RV}, A. Ruiz et L. Vega consid\`erent l'\'equation (\ref{LS}) comme une perturbation de l'\'equation de Schr\"odinger libre pour 
d\'emontrer
(\ref{regularisant}) (en fait une in\'egalit\'e plus forte). Les
potentiels suivants rentrent dans ce cadre:
\begin{itemize}
\item
$\displaystyle V\in L^{p},\; p\geq d/2.$
\item
Des potentiels dans des classes de Morey-Campanato, plus larges que $L^{d/2}$, v\'erifiant une condition de petitesse dans ces espaces (cf \cite{RV} pour les d\'etails). Par exemple, des potentiels admettant des singularit\'es ponctuelles de la forme  
\begin{equation}
\label{potentiel avec pole}
V=\sum_{j=1..N} \frac{a_j}{|x-p_j|^2}+ V_0,\; V_0 \in L^{p},\; p\geq d/2,
\end{equation}
o\`u les $|a_j|$ sont petits.
\end{itemize}
Les potentiels \'etudi\'es dans cet article peuvent aussi d\'ependre
du temps mais nous ne d\'evelopperons pas cet aspect ici. Dans \cite{TD}, (voir aussi \cite{BPSTZ}), l'auteur \'etudie des potentiels de type (\ref{potentiel
avec pole}), (avec $V_0$ r\'egulier) en faisant simplement l'hypoth\`ese de positivit\'e sur
les $a_j$:
$$ a_j+(d/2-1)^2>0,$$
et d\'emontre des in\'egalit\'es sur la r\'esolvante tronqu\'ee:
\begin{equation}
\label{resolvante}
 \forall \lambda \in \CC,\; {\rm Re}\, \lambda^2 \geq C,\; {\rm Im}\,\lambda^2
  \neq 0,\quad ||\chi (P-\lambda^2)^{-1} \chi ||_{L^2\rightarrow L^2} \leq \frac{C}{|\lambda|},\quad
\chi  \in C_0^{\infty}(\RR^d),
\end{equation}
On peut notamment en d\'eduire des estimations de type
(\ref{regularisant}) (cf \cite{BGT}).\par
Les potentiels $V\in L^p$, $p<d/2$, et ceux de la
forme de (\ref{potentiel avec pole}), mais avec des p\^oles d'ordre
sup\'erieur ne sont pas trait\'es dans les travaux pr\'ecit\'es. On
montre ici le caract\`ere critique 
de l'espace $L^{d/2}$ et des singularit\'es en $1/r^2$, $r=|x|$, en construisant un potentiel unipolaire $V$ dont le p\^ole est de l'ordre de $|\log r|^2/r^2$ pr\`es de $0$, tel que $P$ ne v\'erifie ni (\ref{regularisant}) ni (\ref{resolvante}).

\begin{theo}
Soient $d\geq 1$, $N\geq0$ des entiers. Il existe:
\begin{itemize}
\item un potentiel $V$ sur $\RR^d$, positif, radial, \`a support compact tel que:
\begin{gather}
\label{HV1}
V \in C^\infty(\RR^d \backslash \{0\})\\
\label{HV2}
\frac{|\log r|^2}{Cr^2}\leq V(r) \leq \frac{C |\log r|^2}{r^2},\; r\leq r_0.
\end{gather}
\item une suite $(\lambda_n)_{n\geq n_0}$ croissante de r\'eels strictement positifs, de limite $+\infty$;
\item une suite de fonctions $(u_n)_{n\geq n_0}$, radiales, $C^{\infty}$,
  dont le support est de la forme:
\begin{equation*}
\left\{ c_1 \frac{\log \lambda_n}{\lambda_n}\leq r \leq c_2\frac{\log
  \lambda_n}{\lambda_n} \right\}; 
\end{equation*}
\end{itemize}
tels que:
\begin{gather}
\label{equation} 
(-\Delta+V)u_n-\lambda_n^2 u_n=f_n \\
\label{thm norme u}
||u_n||_{L^1}=1\\
\label{thm norme f}
\forall j \in \NN,\;||\frac{d^j}{dr^j}
f_n||_{L^\infty}=O(\lambda_n^{j-N}),\; n\rightarrow +\infty.
\end{gather}
\end{theo}

\begin{corol}
Soient $P=-\Delta+V$, o\`u $V$ est le potentiel
introduit dans le th\'eor\`eme pour un $N\geq 2$, et $\chi\in C_0^{\infty}(\RR^d)$ non nulle en $0$. Alors:
\begin{itemize}
\item (absence d'in\'egalit\'es de Strichartz locales) 
\begin{equation}
\label{reg 1}
\forall q,q_0 \in [1,+\infty],\; q>q_0,\;\forall T>0,\; \forall C>0,\;  
\exists U_0\in C^{\infty}_0,\;||\chi U(t)||_{L^1(]0,T[,L^q(\RR^d))} > C||U_0||_{L^{q_0}(\RR^d)};
\end{equation}
\item (absence d'effet r\'egularisant local)
\begin{equation}
\label{reg 2}
 \forall \sigma>0\;\forall T>0,\; \forall C>0,\;
\exists U_0\in C^{\infty}_0,\; ||\chi U(t)||_{L^1(]0,T[,H^\sigma(\RR^d))} >C ||U_0||_{L^2(\RR^d)};
\end{equation}
\item (absence de dispersion locale)
\begin{equation}
\label{dispersion}
\forall q\in ]1,+\infty],\; \forall T>0,\; \forall C>0,\; \exists
U_0\in C^{\infty}_0,\;
||\chi U(T)||_{L^q(\RR^d)} > C||U_0||_{L^{q'}(\RR^d)};
\end{equation}
\item (absence d'in\'egalit\'es de Strichartz locales avec perte de d\'eriv\'ee)
\begin{multline}
\label{Strichartz}
\forall q\in[1,+\infty],\;\forall \sigma\in[0,1],\; \frac 12 -\frac 1q
>\frac{\sigma}{d},\\ \;\forall T>0,\;\forall C>0,\; \exists U_0\in C^{\infty}_0,\;||\chi U(t)||_{L^1(]0,T[,L^q(\RR^d))}> C||U_0||_{D(P^{\sigma/2})}. 
\end{multline}
\end{itemize}
Dans les in\'egalit\'es pr\'ec\'edentes, on a not\'e $U(t)$ la solution de l'\'equation (\ref{LS}) de condition initiale $U_0$ et $q'$ l'exposant conjugu\'e de $q$ d\'efini par: $1/q+1/q'=1$. 
\end{corol}
\begin{remarques}
\begin{itemize}
\item Si $d>2$, les hypoth\`eses (\ref{HV1},\ref{HV2}) impliquent:
$$ V\in \bigcap_{p< d/2} L^p. $$
\item La suite $(f_n)_{n\geq n_0}$ infirme \'evidemment
(\ref{resolvante}) d\'es que $N\geq 2$.
\item Il appara\^{\i}tra clairement dans la d\'emonstration que l'on peut
construire des quasi-modes d'ordre infini (c'est \`a dire en rempla\c
cant $O(\lambda_n^{j-N})$ par $O(\lambda_n^{-\infty})$ dans (\ref{thm norme
f})) en
s'autorisant une singularit\'e un peu plus forte:
$$ \frac{|\log r|^{2+\eps}}{Cr^2}\leq V(r) \leq \frac{C |\log
r|^{2+\eps}}{r^2},\; \eps>0.$$
En s'autorisant une singularit\'e encore plus grande, on peut imposer \`a la suite $(f_n)$ une d\'ecroissance exponentielle en $-\lambda_n$.
\item Les in\'egalit\'es de Strichartz:
\begin{gather*}
||u||_{L^p\left(]0,T[,L^q(\RR^d)\right)}\leq C ||u_0||_{L^2(\RR^d)},\;
p>2,\; \frac 2p+ \frac dq=\frac d2,
\end{gather*}
sont ni\'ees par (\ref{reg 1}), sauf le cas $q=2$, $p=+\infty$, trivial par
conservation de la norme $L^2$. De m\^eme, (\ref{Strichartz}) nie les
in\'egalit\'es de Strichartz avec perte de d\'eriv\'ee (cf
\cite{BGT}):
\begin{equation*}
||u||_{L^p\left(]0,T[,L^q(\RR^d)\right)}\leq C
||u_0||_{D(P^{\sigma/2})},\;\frac 12-\frac 1q>\frac \sigma d,\;
\sigma\in [0,1],
\end{equation*}
avec des conditions suppl\'ementaires sur $p$, $q$ et $\sigma$. Le cas
limite $\sigma/d=1/2-1/q$ de cette derni\`ere in\'egalit\'e est toujours
vraie par l'inclusion de Sobolev:
$$ L^q(\RR^d) \subset H^{\sigma}(\RR^d)\subset D(P^{\sigma/2}),$$
et la conservation de la norme $D(P^{\sigma/2})$ de toute solution de $(\ref{LS})$.\par
 Remarquons que d'apr\`es cette inclusion, (\ref{Strichartz}) est encore vrai dans les espaces de Sobolev usuels, c'est \`a dire en rempla\c cant $D(P^{\sigma/2})$ par $H^{\sigma}(\RR^n)$. 
\item Le potentiel $V$ et les quasi-modes $u_n$ \'etant \`a supports
compacts, aussi proches de $0$ que l'on veut, le contre-exemple
pr\'ec\'edent est \'egalement valable dans n'importe quel ouvert
born\'e, pour un Laplacien de Neumann ou de Dirichlet.
\item Des exemples d'\'equations de Schr\"odinger ne v\'erifiant pas de
  dispersion ou d'in\'egalit\'es de Strichartz sont \'egalement construits
  dans \cite{Ban} et \cite{CZ}. Dans les deux cas, le probl\`eme vient de la m\'etrique d\'efinissant le Laplacien. Mentionnons aussi que l'id\'ee d'utiliser des quasi-modes en relation avec des in\'egalit\'es de Strichartz est due \`a M. Zworski (voir \cite{SS}). 
\end{itemize}
\end{remarques}
Pr\'ecisons encore que pour le potentiel $V$ du th\'eor\`eme, on peut
d\'efinir $P$ sans ambiguit\'e. L'op\'erateur $-\Delta+V$ a un sens 
sur $C^{\infty}_0(\RR^d \backslash\{0\})$, et est
essentiellement auto-adjoint sur cet espace (cf par exemple le
th\'eor\`eme de Kalf, Walter, Schmincke et Simon cit\'e dans \cite{RS} page
186). On appelle $P$ son unique extension
auto-adjointe, qui est bien s\^ur positive.\par
L'auteur tient \`a remercier Nicolas Burq, son directeur de recherche, et Clotilde Fermanian Kammerer pour leur aide pr\'ecieuse.\par
Le reste du texte est consacr\'e aux d\'emonstrations des r\'esultats
pr\'esent\'es dans cette introduction.

\section{D\'emonstration du corollaire.}
D'apr\`es l'in\'egalit\'e de H\"older, si $v$ est une fonction sur $\RR^d$
avec un support \`a volume fini $B$,
\begin{equation}
\label{volume}
||v||_{L^{q_0}}\leq B^{1/q_0-1/q}||v||_{L^q} ,\quad q>q_0.
\end{equation}
Posons $U_n(t)=e^{-i\lambda_n^2 t} u_n$. Alors:
\begin{gather}
\notag
i\partial_t U_n -P U_n  =  -e^{-i\lambda_n^2 t} f_n,\quad
U_{n\restriction t=0} = u_n\\
\label{expression Un}
U_n(t)=e^{-itP}u_n+i e^{-i\lambda_n^2 t} \int_0^t e^{i(s-t)P} f_n ds.
\end{gather}
Supposons que (\ref{reg 1}) est fausse. 
D'apr\`es (\ref{expression Un}):
$$ \int_0^T ||\chi U_n(t)||_{L^q(\RR^d)}dt\leq \int_0^T ||\chi e^{-itP}u_n||_{L^q(\RR^d)} dt +\int_0^T \int_0^t ||\chi e^{i(s-t)P}f_n||_{L^q(\RR^d)} ds dt.$$
Puisque le support de $u_n$ se concentre en $0$, et que $\chi(0)\neq 0$, le terme de gauche de cette in\'egalit\'e est minor\'e pour $n$ assez grand par:
\begin{equation*}
\int_0^T ||u_n||_{L^q(\RR^d)} dt\geq \frac 1C \left( \frac{\log \lambda_n}{\lambda_n}\right)^{d(1/q-1/q_0)}||u_n||_{L^{q_0}(\RR^d)}. 
\end{equation*}
Quant au terme de droite, il est domin\'e d'apr\`es la n\'egation de (\ref{reg 1}) par:
\begin{equation*}
||u_n||_{L^{q_0}(\RR^d)}+||f_n||_{L^{q_0}(\RR^d)}.
\end{equation*}
On en d\'eduit, en utilisant (\ref{thm norme f}):
$$\left( \frac{\lambda_n}{\log \lambda_n} \right)^{d(1/q_0-1/q)}
||u_n||_{L^{q_0}}=O\left(||u_n||_{L^{q_0}}+\lambda_n^{-N}\right),$$
ce qui donne la contradiction annonc\'ee d\'es que $N\geq 1$ et $q_0<q$, car
d'apr\`es (\ref{thm norme u}), la norme de $u_n$ dans $L^{q_0}$ est
minor\'ee par $1$.\par
Par les in\'egalit\'es de Sobolev, (\ref{reg 1}) implique (\ref{reg
2}) en prenant, dans (\ref{reg 1}), $q_0=2$ et $q>2$ assez proche de 2. Donc (\ref{reg 2})
est vraie aussi.\par
Supposons maintenant la n\'egation de (\ref{dispersion}). Alors, par
(\ref{expression Un}):
\begin{align*}
||\chi u_n(T)||_{L^q(\RR^d)}&\leq C||u_n||_{L^{q'}(\RR^d)} +\int_0^T
||\chi e^{i(s-T)} f_n||_{L^1(\RR^d)} ds\\
&\leq C||u_n||_{L^{q'}(\RR^d)}+C ||f_n||_{L^2(\RR^d)},
\end{align*}
o\`u on a simplement major\'e, pour obtenir la deuxi\`eme
in\'egalit\'e, la norme $L^{q'}$ sur le support de $\chi$ par la norme
$L^2$, \`a une  constante multiplicative pr\`es. On obtient \`a
nouveau une contradiction par (\ref{volume}) et en utilisant $q'<q$.\par
Il reste \`a d\'emontrer (\ref{Strichartz}). Dans le cas contraire, on
aurait, par (\ref{expression Un}):
\begin{equation}
\label{Strichartz un}
\int_0^T ||u_n(t)||_{L^q(\RR^d)} dt \leq
C||u_n||_{D(P^{\sigma/2})}+\int_0^T \int_0^t||e^{i(s-t)P}
f_n||_{D(P^{\sigma/2})} ds\,dt.
\end{equation}
Bien s\^ur, $e^{i(s-t)P}$ envoie isom\'etriquement le domaine de $P^{\sigma/2}$ sur
lui-m\^eme. On a:
$$ ||f_n||^2_{D(P^{1/2})}=\int |\nabla f_n(x)|^2dx+\int (1+V(x))|f_n(x)|^2 dx.$$
Sur le support de $u_n$, donc sur celui de $f_n$, $|x|\geq \frac {\log
\lambda_n}{C \lambda_n}$. En utilisant la majoration de $V$ donn\'ee
par (\ref{HV2}), on obtient que  sur le support de $f_n$:
$$ |V(r)| \leq C \frac{\lambda_n^2}{(\log \lambda_n)^2} \left(
\log\log\lambda_n-\log \lambda_n\right)^2\leq C\lambda_n^{2}.$$
On en d\'eduit, en utilisant (\ref{thm norme f}) avec $N\geq 2$:
\begin{equation}
\label{majoration fn dans D}
||f_n||_{D(P^{1/2})} \leq C\left(||\nabla
f_n||_{L^2}+\lambda_n||f_n||_{L^2}\right) \leq C\lambda_n^{-1}.
\end{equation} 
De plus, l'\'equation $P u_n-\lambda_n^2 u_n=f_n$ implique, en utilisant
\`a nouveau les propri\'et\'es (\ref{thm norme u}) et (\ref{thm norme f}) des normes de $u_n$ et $f_n$:
\begin{align*}
||u_n||^2_{D(P^{1/2})} & \leq ||f_n||_{L^2} ||u_n||_{L^2}+\lambda_n^2
||u_n||^2_{L^2}\\
& \leq C\lambda_n^2 ||u_n||_{L^2}^2.
\end{align*}
Par interpolation sur les normes, on obtient, puisque $\sigma$ est compris entre $0$
et $1$:
\begin{equation}
\label{majoration un dans D}
||u_n||_{D(P^{\sigma/2})}\leq C \lambda_n^\sigma ||u_n||_{L^2}.
\end{equation}
D'apr\`es (\ref{majoration fn dans D}), (\ref{majoration un dans
D}) et l'in\'egalit\'e (\ref{Strichartz un}):
\begin{equation*}
||u_n||_{L^q} \leq C\lambda_n^\sigma ||u_n||_{L^2} +o(1).
\end{equation*}
Donc  avec (\ref{volume}):
\begin{equation*}
\left(\frac{\lambda_n}{\log \lambda_n}\right)
^{d/2-d/q} ||u_n||_{L^2}  \leq C \lambda_n^\sigma ||u_n||_{L^2} +o(1),
\end{equation*}
in\'egalit\'e absurde car $\frac d2-\frac dq>\sigma$.

\section{D\'emonstration du th\'eor\`eme.}

Notons $r=|x|$ et $'$ la d\'eriv\'ee $\frac{d}{dr}$. On cherche $V$, $\lambda_n$, $u_n$ et $f_n$, radiales tels que: 
\begin{equation}
\label{f_n}
 f_n(r)=-u_n''(r)-\frac{d-1}{r} u_n'(r)+V(r) u_n(r) -\lambda_n ^2 u_n(r),
\end{equation}
avec $f_n$ petit \`a l'infini. On commence par changer de fonctions pour \'eliminer la d\'eriv\'ee premi\`ere dans cette expression. En posant:
\begin{equation}
\label{def vn gn W}
u_n=r^{-\frac{d-1}{2}}v_n ,\quad f_n=r^{-\frac{d-1}{2}}g_n,\quad W=V+ \frac{d^2-4d+3}{4 r^2}.
\end{equation}
L'\'equation (\ref{f_n}) devient:
\begin{equation}
\label{g_n}
\tag{\ref{f_n}'}
g_n(r)=-v_n''+W v_n - \lambda_n^2 v_n.
\end{equation}
Soient:
$$y(s)=e^{-\sqrt{s^2+1}},\;
b(s)=-\frac{1}{(s^2+1)^{3/2}}+\frac{s^2}{s^2+1},$$
 solutions $C^{\infty}$ sur $\RR$ de:
\begin{gather}
\label{eqy}
-y''(s)+b(s)y(s)=0\\
\label{inegy}
y(s)>0,\forall j\in \NN,\quad |y^{(j)}(s)|\leq C_je^{-|s|}\\
\label{inegb}
|b(s)|\leq 1.
\end{gather}
On notera, pour $2$ param\`etres r\'eels $\omega$ et $a$:
\begin{gather*}
y_{\omega,a}(r)=y(\omega(r-a)),\quad b_{\omega,a}(r)=\omega^2 b(\omega(r-a)),
\end{gather*}
solutions de:
\begin{equation}
\label{eqy'}
\tag{\ref{eqy}'}
-y''_{\omega,a}+b_{\omega,a} y_{\omega,a}=0;
\end{equation}
Soit $q(\lambda)>0$ d\'efinie au voisinage de $+\infty$,
  strictement croissante et v\'erifiant:
\begin{gather}
\label{hypq1}
\lim_{\lambda \rightarrow +\infty} q(\lambda)=+\infty\\
\label{hypq2}
\lim_{\lambda \rightarrow +\infty} \frac{q(\lambda)}{\lambda}=0.
\end{gather}
On d\'efinit la suite $(\lambda_n)_{n\geq n_0}$, qui tend vers $+\infty$ plus vite que $10^n$ par:
$$ 10^n=q(\lambda_n).$$
Soit une troncature:
$
\displaystyle
\chi \in C^{\infty}_0\left(\left]-1,1\right[ 
\right),\; \chi=1 \;{\rm sur }\; \left[-\frac 12,\frac 12\right].
$\par
On notera:
\begin{gather*}
\psi_n(r)=\chi(10^n r)-\chi(10^{n+1}r),\\
\chi_n(r)=\chi\left( 10^{n+1}\left(r-\frac{3}{10^{n+1}}\right)\right),
\end{gather*}
\setlength{\unitlength}{1mm}
\begin{center}
\begin{figure}
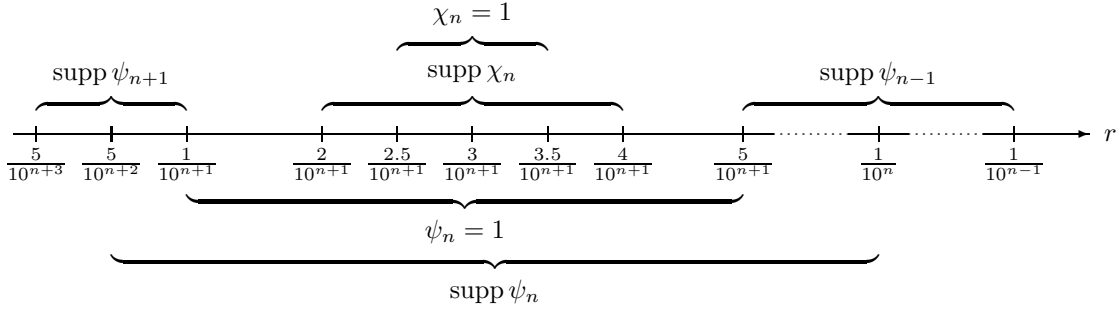

\dessineloc
\caption{Troncatures autour de $r=3/10^n+1$}
\end{figure}
\end{center}
Les $\psi_n$ forment une partition de l'unit\'e au voisinage de $0$, et chacun des $\chi_n$ localise sur le support de $\psi_n$, en dehors du support des autres $\psi_j$, $j\neq n$:
\begin{gather}
\notag
\supp \psi_n\subset \left[\frac{5}{10^{n+2}},\frac{1}{10^{n}}\right]\\
\notag
\frac {1}{10^{n+1}}\leq r\leq\frac {5}{10^{n+1}} \Rightarrow \psi_n(r)=1\\
\label{somme psin}
\sum_{n\geq n_0} \psi_n(r)=\chi(10^{n_0}r)\\
\label{support chin}
\supp \chi_n \subset \left[ \frac{2}{10^{n+1}},\frac {4}{10^{n+1}} \right]\subset \{\psi_n=1\}\\
\label{chin=1}
\frac {25}{10^{n+2}}\leq r\leq \frac {35}{10^{n+2}} \Rightarrow \chi_n(r)=1.
\end{gather}
On notera $y_n$ et $b_n$ les fonctions:
\begin{gather}
\label{def yn}
 y_n(r)=y_{\frac{\lambda_n}{2},\frac{3}{10^{n+1}}}(r)=y\left( \frac {\lambda_n}{2}\left(r-\frac{3}{10^{n+1}}\right)\right)\\
 b_n(r)=b_{\frac{\lambda_n}{2},\frac{3}{10^{n+1}}}(r)=\frac{\lambda_n^2}{4} b\left(\frac{\lambda_n}{2}\left(r-\frac{3}{10^{n+1}}\right)\right).
\end{gather}
On peut voir chacun des $-b_n$ comme un puits de potentiel qui, par l'\'equation (\ref{eqy'}), concentre l'\'energie de $y_n$ autour du point $3/10^{n+1}$. On va construire $W$ pour que l'\'equation $g_n(r)=0$ s'\'ecrive, sur  un petit intervalle autour de $3/10^{n+1}$ (qui inclura le support de $v_n$) $-v_n''+b_n v_n=0$. La taille de cette intervalle sera de l'ordre de $10^{-n}$, tendant moins vite vers $0$ que le param\`etre $\omega=\lambda_n/2$ de l'\'equation. En d\'efinissant $v_n$ comme une troncature de $y_n$, la
d\'ecroissance exponentielle de $y$ \`a l'infini induira gr\^ace \`a la
diff\'erence d'ordre de ces deux \'echelles, des termes d'erreur
assez petits pour impliquer (\ref{thm norme u}) et (\ref{thm norme f}).
\begin{center}
\begin{figure}
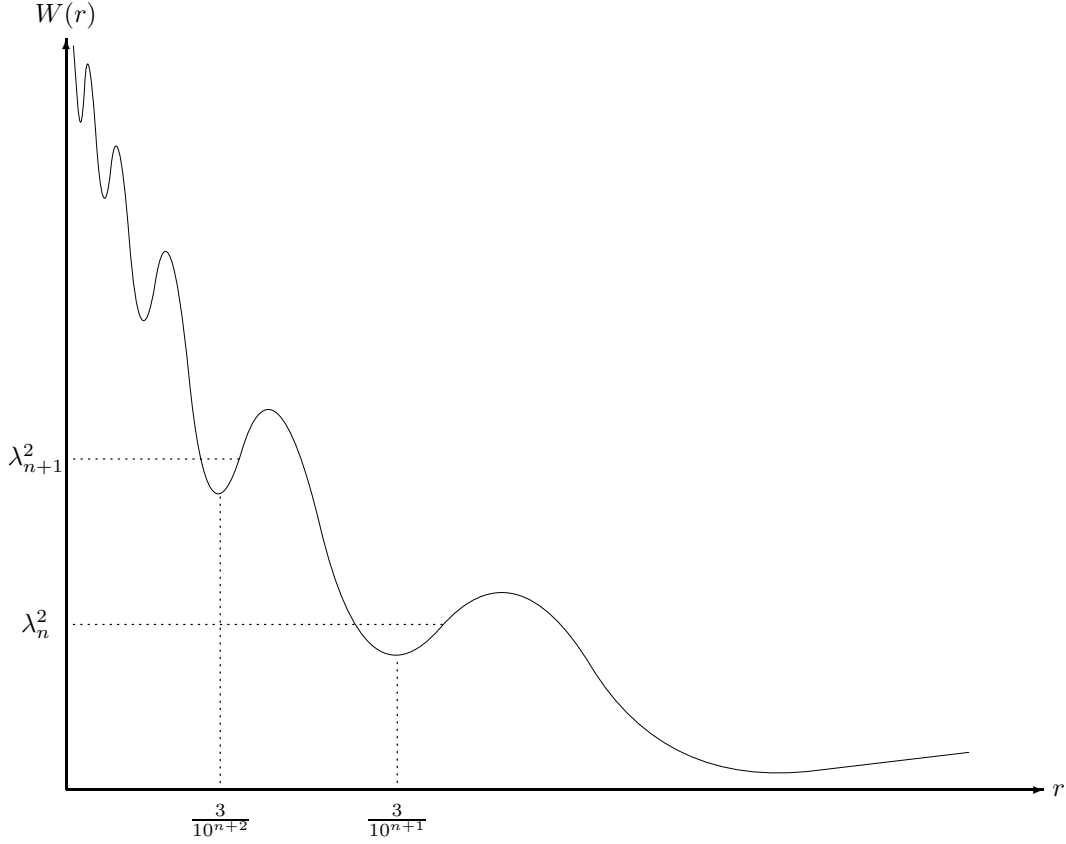

\label{figure 1}
\dessineW

\caption{Le potentiel W}
\end{figure}
\end{center}
Posons (cf figure):
\begin{gather}
\label{def W}
W(r)=\sum_{n\geq n_0} \psi_n(r)
\left(b_n(r)+\lambda_n^2\right)\\
v_n(r)= \alpha_n y_n(r)\chi_n(r),
\end{gather}
o\`u les $\alpha_n$ sont des constantes strictement positives \`a choisir. Les fonctions $g_n$, $u_n$ et $f_n$, ainsi que le potentiel $V$ sont alors d\'efinies par (\ref{g_n}) et (\ref{def vn gn W}). 
\begin{lem}
\label{estimations}
On a les in\'egalit\'es suivantes:
\begin{gather}
\label{ineg fn}
\forall j\in \NN,\;\exists C>0,\; ||\frac{d^j f_n}{dr^j}||_{L^{\infty}(\RR^d)}\leq C \alpha_n \lambda_n^{1+j}\left(q(\lambda_n)\right)^{\frac{d+1}{2}} 
 e^{-\frac{\lambda_n}{20 q(\lambda_n)}}\\
\label{ineg un}
||u_n||_{L^1(\RR^d) } \geq \frac{\alpha_n}{C} \lambda_n^{-1} \left(q(\lambda_n)\right)^{-\frac{d-1}{2}}.
\end{gather}
\end{lem}
\begin{proof}
D'apr\`es (\ref{support chin}), sur le support de $\chi_n$, $W$ vaut exactement $b_n+\lambda_n^2$. On en d\'eduit:
\begin{equation}
\label{expression gn}
g_n(r) =-v_n''(r) + b_n(r) v_n(r).
\end{equation}
Soit, en utilisant l'\'equation (\ref{eqy'}):
\begin{equation*}
f_n(r)=r^{-\frac{d-1}{2}} g_n(r)=r^{-\frac{d-1}{2}}\alpha_n \left(y_n(r) \chi_n''(r) +2 y'_n(r) \chi_n'(
r)\right).
\end{equation*}
La d\'eriv\'ee $j$-\`eme de $f_n$ est donc de la forme:
\begin{equation}
\label{derivee fn}
\frac{d^j f_n}{dr^j}=\sum_{j_1+j_2+j_3=j+1}
\beta_{j_1,j_2,j_3}\frac{d^{j_1}}{dr^{j_1}}y_n 
\frac{d^{j_2}}{dr^{j_2}}\chi_n'\frac{d^{j_3}}{dr^{j_3}} r^{-\frac{d-1}{2}}. 
\end{equation}
Sur le support de $\chi'_n$, on a:
$$ \left|r- \frac{3}{10^{n+1}}\right| >
\frac{1}{10^{n+1}}=\frac{1}{10q(\lambda_n)}.$$ 
Donc d'apr\`es les majorations (\ref{inegy}):
\begin{align*}
|y_n^{(j_1)}(r)|&\leq C\lambda_n^{(j_1)} e^{-\frac{\lambda_n}{2}\left|r-
\frac{3}{10^{n+1}}\right|}\\
&\leq \lambda_n^{j_1} e^{-\frac{\lambda_n}{20q(\lambda_n)}}
\end{align*}
De plus:
\begin{alignat*}{3}
|\chi_n^{(j_2+1)}|&\leq C\left(q(\lambda_n)\right)^{j_2+1}&\leq& C\lambda_n^{j_2}q(\lambda_n)\\
\left|\left(\frac{d}{dr}\right)^{j_3}\left(
    r^{-\frac{d-1}{2}}\right)\right|&\leq
C\left(q(\lambda_n)\right)^{\frac{d-1}{2}+j_3}&\leq& C\lambda_n^{j_3}\left(q(\lambda_n)\right)^{\frac{d-1}{2}},
\end{alignat*}
car sur le support de $\chi_n$, $r\geq \frac{1}{10^{n+1}}$. En majorant la
somme (\ref{derivee fn}) par ces trois derni\`eres in\'egalit\'es, on
obtient exactement (\ref{ineg fn}).
(\ref{derivee fn}), on obtient exactement (\ref{ineg fn}).\par
D\'emontrons (\ref{ineg un}). Par d\'efinition de $y_n$:
$$ \frac{\lambda_n}{2} \left | r-\frac{3}{10^n} \right|\leq \frac 12 \Rightarrow y_n(r)\geq m=\sup_{|s|\leq \frac 12} |y(s)|.$$
De plus, si cette derni\`ere condition sur $r$ est v\'erifi\'ee, et si $n$ est assez grand on a $\chi_n(r)=1$ et donc:
$$ u_n(r)=r^{-\frac{d-1}{2}} v_n(r)=r^{-\frac{d-1}{2}} \alpha_n y_n(r)$$
D'o\`u:
\begin{align*}
||u_n||_{L^1} &\geq m \alpha_n \int_{
\lambda_n \left|r-\frac{3}{10^n} \right|\leq 1} r^{-\frac{d-1}{2}} r^{d-1} dr\\
&\geq \frac{\alpha_n}{C} \left(10^{-n}\right)^{\frac{d-1}{2}} \lambda_n^{-1}. 
\end{align*}
\end{proof}
Choisissons un $M>1$ et posons: 
\begin{equation}
\label{defq}
q(\lambda)= \frac{\lambda}{M\log \lambda}. 
\end{equation}
La fonction $q$, positive et strictement croissante au voisinage de $+\infty$, v\'erifie bien les hypoth\`eses (\ref{hypq1}) et (\ref{hypq2}). De plus, d'apr\`es le lemme \ref{estimations}, et en majorant grossi\`erement $q(\lambda_n)$ par $\lambda_n$:
\begin{gather*}
||\frac{ d^jf_n}{dr^j}||_{L^{\infty}} \leq C \alpha_n \lambda_n^{\frac{d+3}{2}+j-\frac{M}{20}}\\
||u_n||_{L^1} \geq \frac{\alpha_n}{C} \lambda_n^{\frac{d+1}{2}},
\end{gather*}
avec des nouvelles constantes $C$ qui d\'ependent \'eventuellement de $M$.
Les conditions (\ref{thm norme u}) et (\ref{thm norme f}) du th\'eor\`eme
sont donc respect\'ees avec un bon choix des constantes $M$ et
$\alpha_n$. Le support de $u_n$ est confondu avec celui de $\chi_n$, qui
est effectivement de la forme:
$$\left \{ c_1\frac{\log \lambda_n}{\lambda_n}\leq r \leq c_2
  \frac{\log\lambda_n}{\lambda_n}\right\},$$
si $\chi$ est bien choisie.
Il reste \`a v\'erifier l'assertion (\ref{HV2}) sur le potentiel. On
commence par \'evaluer la fonction r\'eciproque de $q$.
\begin{lem}
\label{equiv q}
Soit $q$ est d\'efinie par (\ref{defq}). Alors:
$$ q(\lambda)\log(q(\lambda)) \sim \frac{\lambda}{M},\; \lambda\rightarrow +\infty.$$
\end{lem}
\begin{proof}
\begin{align*}
q(\lambda)\log q(\lambda)&=q(\lambda) \log\left(\frac{\lambda}{M\log \lambda}\right)\\
&=q(\lambda)\left(\log \lambda-\log\log \lambda -\log M\right)\\
&\sim q(\lambda)\log\lambda,
\end{align*}
quand $\lambda$ tend vers $+\infty$.
\end{proof}
Sur le support de $\psi_n$, puisque $q(\lambda_n)=10^n$,
\begin{gather*}
 \frac{1}{20r}\leq q(\lambda_n) \leq \frac{1}{r}\\
\frac{1}{20r}|\log r -\log 20| \leq q(\lambda_n) \log(q(\lambda_n))\leq \frac{1}{r} |\log r|,
\end{gather*}
(l'application $s\mapsto s \log s$ \'etant croissante pour $s>e^{-1}$).
Soit, en utilisant le lemme \ref{equiv q} et en se pla\c cant assez pr\`es de $0$:
$$ C^{-1} \frac 1r |\log r|\leq \lambda_n \leq C \frac{1}{r} |\log r|. $$
Donc, par la d\'efinition (\ref{def W}) de $W(r)$, et comme $|b_n|$ est major\'e par $\lambda_n^2/4$:
$$ C^{-1} \frac {|\log r|^2}{r^2} \sum_{n\geq n_0} \psi_n(r) \leq W(r) \leq C  \frac {|\log r|^2}{r^2} \sum_{n\geq n_0} \psi_n(r), $$
ce qui donne l'in\'egalit\'e (\ref{HV2}) sur le potentiel $W$ compte tenu de (\ref{somme psin}). On en d\'eduit \'evidemment (\ref{HV2}) sur $V$.\par
\begin{remarque} 
Pour avoir des quasi-modes d'ordre infini il suffisait de prendre $q(\lambda)=\frac{\lambda}{(\log\lambda)^{1+\eps}}$, et de modifier le lemme \ref{equiv q} en cons\'equence.
\end{remarque}